\documentclass[final]{siamltex}
\usepackage{bm}
\usepackage{epsfig}
\usepackage{amsfonts}
\usepackage{color}

\begin{document}

\title{Uncertainty Quantification for Atomistic Reaction Models:
        An Equation-Free Stochastic Simulation Algorithm Example}

\author{Yu Zou 
        \and Ioannis G. Kevrekidis \thanks{Department of Chemical Engineering and Program in
Applied and Computational Mathematics, Princeton University, Princeton, NJ 08544({\tt yannis@princeton.edu}).}}

\maketitle

\begin{abstract}

We describe a computational framework linking
Uncertainty Quantification (UQ) methods for continuum problems
depending on random parameters with Equation-Free (EF) methods for
performing continuum deterministic numerics by acting directly on
atomistic/stochastic simulators.
Our illustrative example is a heterogeneous catalytic reaction
mechanism with an uncertain {\em atomistic} kinetic parameter; the
``inner" dynamic simulator of choice is a Gillespie Stochastic
Simulation Algorithm (SSA).
We demonstrate UQ computations at the coarse-grained level
{\em in a nonintrusive way}, through the design of brief, appropriately
initialized computational experiments with the SSA code.
The system is thus observed at three levels: (a) a fine scale for each
stochastic simulation at each value of the uncertain parameter;
(b) an intermediate coarse-grained state
for the {\it expected behavior} of the SSA at each value of the uncertain
parameter; and (c) the desired fully coarse-grained level:
distributions of the coarse-grained
behavior {\em over the range of uncertain parameter values}.
The latter are computed in the form of generalized Polynomial Chaos (gPC)
coefficients in terms of the random parameter.
Coarse projective integration and coarse fixed point computation
are employed to accelerate the computational evolution of these
desired observables, to converge on random stable/unstable steady states,
and to perform parametric studies with respect to other (nonrandom) system
parameters.

\end{abstract}

\begin{keywords}
Uncertainty quantification, Equation-Free,
generalized Polynomial Chaos, stochastic Gillespie algorithm, multiscale
\end{keywords}

\section{Introduction}
The temporal evolution of many engineering systems can be described through continuum models
(typically Ordinary or Partial Differential Equations, ODEs or PDEs) depending on
{\em random} parameters and/or initial/boundary conditions.
A number of methods have been proposed to study the evolution of the probability distribution
of the solutions of such random problems (the so-called uncertainty quantification or UQ).
Among early exploration approaches in this area are Monte-Carlo
based methods \cite{Metropolis:49,Shinozuka:72,Bruns:81,Fishman:96},
which normally require a large number of ensemble realizations to
achieve convergence.
For slightly perturbed systems, whose random parameters can be described as
small fluctuations around average values, perturbation methods may
be utilized \cite{Boyce:64,Hart:70,LiuA:86}.
Such methods cannot, however, be used to treat more generic random systems,
exhibiting large parametric uncertainties, and even if used they can only obtain information
on low-order statistics.
To overcome this difficulty, moment-closure techniques
(e.g., \cite{Bharrucha-Reid:68}) were tried to study dynamics of statistical moments of solutions for which
the small-uncertainty assumption in the perturbation method may be relaxed.
A significant difficulty with this aprroach lies in deriving an
accurate {\em closure} for high-order statistical moments of the
solution distributions; this is extremely difficult, especially for
nonlinear systems.

In recent years an alternative approach for UQ, the stochastic Galerkin method,
has received considerable attention as an approach to solving
ODEs or PDEs with random parameters.
This approach originated from the work of Wiener \cite{Wiener:38}
in constructing multiple stochastic integrals (also known as Homogeneous Chaos)
to represent functionals of Wiener processes.
This idea was then utilized in \cite{Ghanem:91} to express
solutions of random systems in terms of Hermite polynomials of Gaussian random variables
(i.e., Wiener-Hermite Polynomial Chaos).
Projections of these solutions on the Hermite polynomial basis
are deterministic, and can be readily numerically solved after truncation
through a Galerkin method.
Early works in applying the method in engineering systems include \cite{Ghanem:91},
where random static and dynamical problems in structural analysis were investigated.
The method was subsequently applied in uncertainty quantification
of various physical systems including (but certainly not restricted to) porous media \cite{Ghanem:98},
fluid dynamics \cite{Maitre:01} and chemical reactions \cite{Reagan:04}.
The application of the Wiener-Hermite Polynomial Chaos was recently extended in \cite{XiuB:02,XiuA:02,Xiu:03}
to a general situation, where the Askey scheme was used to
construct {\em generalized} Polynomial Chaos (gPC) with respect to
other continuous and discrete probability measures.
Other developments in the application of the method include
piecewise function representation \cite{Deb:01},
non-orthogonal expansions \cite{Babuska:02}
and Wiener-Haar wavelet expansions \cite{Maitre:04}.
One advantage of the stochastic Galerkin method relative to
Monte Carlo simulation is that it can reduce a random system to a deterministic one,
often with considerably fewer degrees of freedom.
Moreover, this method -- along with its extensions and ramifications --
can handle more general uncertainty problems such as large
fluctuations, multimodal and discontinuous probability distributions, etc.

In order to implement the stochastic Galerkin method, we must be able to derive differential equations
for the projection of solution distributions onto gPC bases (i.e., ODEs for the evolution of
gPC coefficients).
When explicit equations for these gPC coefficients are not available in closed form,
a unified approach, combining equation-free techniques \cite{Theodoropoulos:00,Kevrekidis:03,Kevrekidis:04}
and the stochastic Galerkin method {\em in a nonintrusive way} has been proposed.
This is {\em equation-free uncertainty quantification}
\cite{Xiu:05}, which can evolve gPC coefficients in time, and has
been used to perform steady state and limit cycle bifurcation
analysis of the gPC coefficient equations without needing these
equations in closed form.
This approach is based on short bursts of direct simulation (ensembles of such simulations
distributed over the random parameter(s)) as an inner,  ``fine scale" simulator; these bursts
are used to numerically estimate the necessary information at the ``coarser" gPC coefficient level
(e.g., temporal derivatives of gPC coefficients) {\em on demand}.
The main assumption underlying this approach is that
the long-term dynamics in gPC coefficient space lie on a low-dimensional,
attracting manifold, which can be parametrized by only a few leading-order
gPC coeffiicents in the appropriate orthogonal polynomial chaos basis (see \cite{XiuB:02}).
In this way, model reduction can be achieved by using gPC expansions.

When we study continuum models of chemically reacting systems (e.g. ODEs for the evolution of
reactant concentrations in a stirred tank reactor, or for the evolution of coverages on
catalytic surfaces), UQ techniques can be used to study the effect of uncertainty in
kinetic or operating parameters on the overall reactor behavior (e.g. steady state
concentrations, reaction rates etc.).
In many problems of contemporary interest, however, such continuum differential equations
are not available in closed form; instead, we are given a description of the chemical kinetics
at an atomistic/stochastic level and the uncertainty enters now in the parameters of the
atomistic/stochastic simulation itself (e.g. in certain transition probabilities).
The normal steps for uncertainty quantification modeling would involve (a) derivation
of closed continuum kinetic equations for the coarse-grained observables of the
atomistic simulation (averaged concentrations, averaged coverages)-- these equations need
to explicitly express {\em coarse-grained} parameter uncertainty in terms of the
{\em fine scale, atomistic level} parameter uncertainty; and (b) the application
of traditional, continuum UQ techniques on these continuum kinetic equations.
Here we will show how to circumvent this ``intermediate" derivation of coarse-grained,
continuum kinetic equations.
We will implement a computational framework that performs UQ computations for these
unavailable, random continuum kinetic equations {\em by acting on the fine scale,
atomistic/stochastic simulator directly}.
This uncertainty quantification procedure for stochastic chemical
reaction models with uncertain parameters involves three distinct
system levels: a fine scale, microscopic simulator {\em for fixed
values of the uncertain parameter(s)}; an {\em intermediate}
coarse-grained scale, where continuum kinetic ODEs exist {\em in
principle} for certain system observables, still for fixed values of
the uncertain parameter(s); and, a final, fully coarse-grained scale
of gPC coefficients for the distributions of the solutions of these
(unavailable) closed form continuum equations.
The equation-free machinery is thus used at two successive levels: first to circumvent
the derivation of closed form continuum kinetic equations; and then to circumvent
the derivation of explicit equations for the evolution of the gPC coefficients for
the solutions of these equations with uncertain parameters through a Galerkin procedure.
The paper is organized as follows: We start in Section \ref{StoGalerkin:sec} briefly
recalling the formulation of the stochastic Galerkin method.
Equation-free techniques are combined with this stochastic
Galerkin method in Section \ref{EFPIFP:sec} to study random systems {\em without explicitly
available governing equations}.
The approach is illustrated through the standard model for the
$A+1/2B_2 \rightarrow AB$ reaction, which has been used as a
simplified description for the catalytic oxidation of CO
\cite{Makeev:02}.
The fine-scale simulator is chosen to be the
Gillespie Stochastic Simulation Algorithm (SSA).
Computational results are presented in Section \ref{results:sec}
where we also demonstrate the efficiency of using the Gauss-Legendre
quadrature for computing gPC projections.
Random stable/unstable steady-state solutions are also computed
in Section \ref{results:sec}, and continuation algorithms are
implemented to explore the effect of variation of a system
(nonrandom) parameter on the statistics of the random steady state.
We conclude with a brief summary and discussion in Section \ref{conclusion:sec}.
\section{The Stochastic Galerkin Method} \label{StoGalerkin:sec}
Consider a system whose evolution in time is governed by
 the differential equation
\begin{equation}
  { {d{\bm x}_c} \over {dt} } = {\bm f} ({\bm x}, \omega_c), \quad {\bm x}_c(\omega_c,0) = {\bm x}_{c,0}(\omega_c);
\label{ODE:eqn}
\end{equation}
the system state is ${\bm x}_c(\omega_c,t)$, where $\omega_c$ is an
element in the sampling space $\Omega_c$.
The subscript $c$ stands for ``coarse", since we consider this to be the coarse-grained
description of an atomistic simulator (an SSA simulator in our illustrative example below).
The solution of the above equation can be described by an
expansion in an $L^2$ space with a generalized polynomial chaos basis $\Psi_i({\bm \xi}(\omega_c))$, i.e.,
\begin{equation}
   {\bm x}_c(\omega_c,t) = \sum_{i=0}^P {\bm x}_{cc}^i(t) \Psi_i({\bm \xi}(\omega_c)).
\label{PCexpansion:eqn}
\end{equation}
The projections ${\bm x}_{cc}^i(t)$ are determined by
\begin{equation}
   {\bm x}_{cc}^i(t) = { {<{\bm x}_c(\omega_c,t),\Psi_i({\bm \xi}(\omega_c))>} \over {<\Psi_i({\bm \xi}(\omega_c)),\Psi_i({\bm \xi}(\omega_c))>}},
\label{PCcoefficient:eqn}
\end{equation}
where the inner product $<\cdot, \cdot>$ is defined as
\begin{equation}
  <q({\bm x}({\bm \xi})), g({\bm x}({\bm \xi}))> = \int q({\bm x}({\bm \xi})) g({\bm x}({\bm \xi})) p({\bm \xi}) d {\bm \xi}, \nonumber
\end{equation}
(where $p({\bm \xi})$ is a probability measure of ${\bm \xi}$) for any
functions $q({\bm x}({\bm \xi}))$ and $g({\bm x}({\bm \xi}))$ in the $L^2$ space.
The subscript ``cc" is indicative of a second level of coarse-graining: these coefficients describe
the {\em statistics} of the distribution of solutions for the already coarse-grained state ${\bm x}_c(\omega_c,t)$
over the sampling space.

Through a Galerkin projection, equations governing gPC coefficients
${\bm x}_{cc}^i(t)$ are obtained as
\begin{eqnarray}
 &&  { {d {\bm x}_{cc}^i} \over {dt} } = { {< {\bm f} (\sum_{i=0}^P {\bm x}_{cc}^i(t) \Psi_i({\bm \xi})), \Psi_i({\bm \xi})>} \over {<\Psi_i({\bm \xi}),\Psi_i({\bm \xi})>} } , \nonumber
      \\
 &&  i=0,2,\cdots,P,
\label{coarserODE:eqn}
\end{eqnarray}
with ${\bm x}_{cc}^i(0) =   { {<{\bm x}_c(0),\Psi_i({\bm \xi})>} \over {<\Psi_i({\bm \xi}),\Psi_i({\bm \xi})>}}$.
The above equation can be formally compactly rewritten as
\begin{equation}
   { {d {\bm X}_{cc}} \over {dt} } = {\bm H}( {\bm X}_{cc} ),
\label{Galerkin:eqn}
\end{equation}
where
\begin{equation}
  {\bm X}_{cc}=({\bm x}_{cc}^0, {\bm x}_{cc}^1, \cdots,{\bm x}_{cc}^P)^T \nonumber
\end{equation}
and
\begin{equation}
 {\bm H}=({\bm h}_0, {\bm h}_1, \cdots, {\bm h}_P)^T. \nonumber
\end{equation}
Here
\begin{eqnarray}
  && {\bm h}_i({\bm X}_{cc})= { {< {\bm f} (\sum_{i=0}^P {\bm x}_{cc}^i(t) \Psi_i({\bm \xi})), \Psi_i({\bm \xi})>} \over {<\Psi_i({\bm \xi}),\Psi_i({\bm \xi})>} } , \nonumber
   \\
  && i=0,2,\cdots,P.
\end{eqnarray}
If ${{d {\bm X}_{cc}} \over {dt}} = {\bm 0}$ in the long-time limit,
then Equ.(\ref{Galerkin:eqn}) has a steady state, which can be used
to obtain the probability distribution of the random steady state of
(\ref{ODE:eqn}).
The explicit derivation of Equation (\ref{Galerkin:eqn}) is a challenging problem
if (\ref{ODE:eqn}) is a set of strongly nonlinear equations; pseudospectral approaches (see e.g. \cite{Reagan:04})
provide a possible alternative.

\section{Equation-Free Computation for Random Dynamical Systems without Explicit Governing Equations} \label{EFPIFP:sec}
Equation-free methods have been applied in recent years to investigate
solutions of {\em non-random} macroscopic systems whose evolution equations are not
explicitly available \cite{Theodoropoulos:00,Makeev:02,Kevrekidis:03,Kevrekidis:04,Zou:05}.
The approach can, in principle, provide clear scenarios of the
coarse-level evolution and its parametric dependence while requiring
only short-time bursts of evolution with the micro-level simulators;
in effect, it is a framework of accelerating the extraction of
information from the microscale simulation through judicious design
of computational experiments and processing of their results.

The equation-free approach utilizes the so-called coarse time-stepper as its basic element; this
time-stepper consists essentially of three components: {\it lifting},
{\it micro-simulation}, and {\it restriction}.
Lifting is a protocol that transforms a coarse-level state to consistent fine-level states;
restriction is the converse of lifting.
Note that the lifting will, in general, not be a one-to-one transformation,
since fine-scale states have far more degrees of freedom than their
corresponding coarse-grained descriptions; this is a vital step in EF
computations, and a good code should test, on line, that different realizations
of the lifting protocol do not affect the coarse-grained computational results.
This approach was applied recently in conjunction with the stochastic Galerkin method to
study random dynamical systems \cite{Xiu:05}.
The fine scale state, in this case, is a (large enough) ensemble of system realizations;
at each time snapshot, this ensemble can be represented by its projection
onto a generalized polynomial basis; if a low order gPC truncation provides an accurate representation,
these first few gPC coefficients constitute the coarse-grained description.
In principle, there exist differential equations for these projections,
which can be derived, implemented in a computer code, and then solved using numerical methods.
However, it may be difficult or even impossible to derive these equations;
equation-free methods were utilized as an alternative way of solving them without deriving them first.
Assuming that the coarse-grained evolution is smooth,
one can use techniques like projective integration \cite{GearA:03,GearB:03} to accelerate the
time-evolution of the gPC coefficients using fine scale (ensemble realization)
simulations over only relatively short time segments.
This is accomplished by {\em observing} the evolution of the
ensemble MC runs on their gPC coefficients (through {\em
restriction}), and use of these data to locally estimate the
time-derivative of the coarse-grained description (the local time
derivative of the gPC coefficients).
This information is then ``passed" to traditional continuum numerical
initial value problem solvers (ranging from the simple explicit forward
Euler to Runge-Kutta type or even implicit integrators) that ``project"
the coarse-grained state forward in time \cite{Kevrekidis:03,Kevrekidis:04}.
One thus solves the initial value problem for the coarse-grained description
with the necessary quantities (the gPC local time derivatives) obtained
not through a function evaluation from a closed-form model, but through
processing the results of ``judicious" numerical experiments with the fine
scale code.
Beyond coarse projective integration,
when Equation (\ref{coarserODE:eqn}) has a stationary state,
one can turn the coarse time-stepper into a fixed-point operator
and use matrix-free methods such as Newton-GMRES \cite{Kelley:95}
to compute  stable/unstable coarse-grained steady states.
These correspond to random steady states of the original random equations.

The equation-free technique in \cite{Xiu:05}, however, requires the evolution
equations for the random dynamical systems to be explicitly available.
In the case that these equations are not explicitly available,
we now show how to exploit fine-scale models {\em underlying} these differential equations.
We now have two successive lifting (and, correspondingly, two
successive restriction) levels.
To obtain a numerical representation for evolution of the desired, ``doubly" coarse-grained
representation (the gPC coefficients), the fine-scale states of these
fine-scale models are first restricted to an intermediate coarse level: we obtain
individual states in the ensemble of random ODE solutions.
At a second level of restriction, the entire ensemble of these states is used
to compute their gPC coeffcients.
We view the level of random differential equation as the ``intermediate coarse"
scale and the level where the gPC projections reside as a desired ``fully coarse" scale.

The interaction between intermediate coarse and the fully coarse scale states
is embodied in equations (\ref{PCexpansion:eqn}) and (\ref{PCcoefficient:eqn}),
which correspond to {\it lifting} and {\it restriction}, respectively,
between these two scales.
The interaction between fine and ``intermediately coarse" states {\em at fixed values
of the random parameters $\bm{\xi}$} is described by
two new lifting and restriction mappings: $\bm \mu$ ({\it lifting}) and $\cal M$ ({\it restriction}),
\begin{equation}
   {\bm X}_f(\bm{\xi},t) = {\bm \mu} ( {\bm x}_c(\bm{\xi},t)),
\label{lifting:eqn}
\end{equation}
\begin{equation}
   {\bm x}_c(\bm{\xi},t) = {\cal M} ({\bm X}_f(\bm{\xi},t)).
\label{restriction:eqn}
\end{equation}
Recall that
the fine-scale states ${\bm X}_f$ are characterized by many more degrees of freedom
than their intermediate-scale   ${\bm x}_c$ counterparts.
As the fine-scale states are restricted to the intermediate-scale level, these additional
degrees of freedom are eliminated.

\begin{figure*}
\centerline{\epsfig{figure=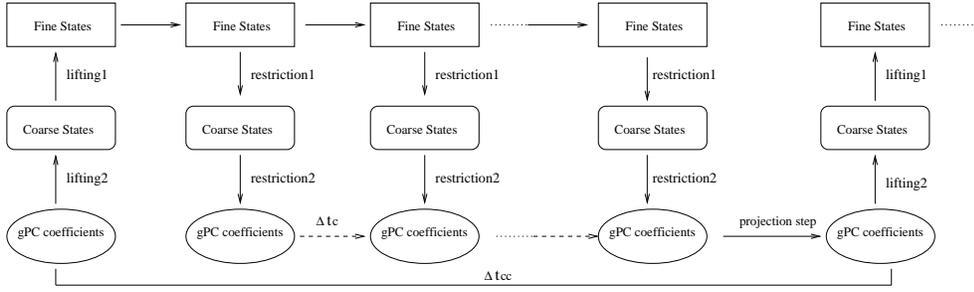,width=\textwidth}}
\caption{Schematic of coarse projective integration for  multiscale
systems with uncertainty.} \label{PImultiscale:fig}
\end{figure*}

Coarse projective integration for random systems whose coarse-scale equations are
not explicitly available can be summarized in the following steps: (see Fig. \ref{PImultiscale:fig})

\begin{enumerate}
\item Generate an ensemble of (intermediate) coarse-scale states based on initial values of gPC
    coefficients by (\ref{PCexpansion:eqn}).
\item Generate an ensemble of corresponding fine-scale states consistent
    with each element of the ensemble of (intermediate) coarse-scale states
    by (\ref{lifting:eqn}). In the case that (intermediate) coarse-scale states are mean fields,
    the average of these fine-scale states should be equal to the prescribed (intermediate) coarse state.
\item Evolve fine-scale states via fine-scale (in the example below, SSA) simulators.
\item Restrict the fine-scale states to the (intermediate) coarse level and obtain
    an ensemble of these coarse-scale states by (\ref{restriction:eqn}).
\item Further restrict (intermediate) coarse-scale states to the desired fully coarse level,
    to obtain gPC coefficients by (\ref{PCcoefficient:eqn}).
\item Perform the above two steps successively, and use the results to estimate the temporal
    derivative of the fully coarse observables (gPC coefficients). 
    The data collection time is dictated by the separation of fast/slow time scales we assume 
     prevails at the coarse grained level; also by the noise of gPC coefficients brought in by 
     the SSA at the fine level and the approach used to compute these coefficients at the coarse level.
\item {\em Project forward in time} --using a continuum numerical integrator, such as forward Euler-
     to obtain the fully coarse observables at a later time. Go back to Step 1. The selection of the projective time step so as to retain overall stability of the projective method is discussed in \cite{GearA:03}.
\end{enumerate}

\begin{figure}
\centerline{\epsfig{figure=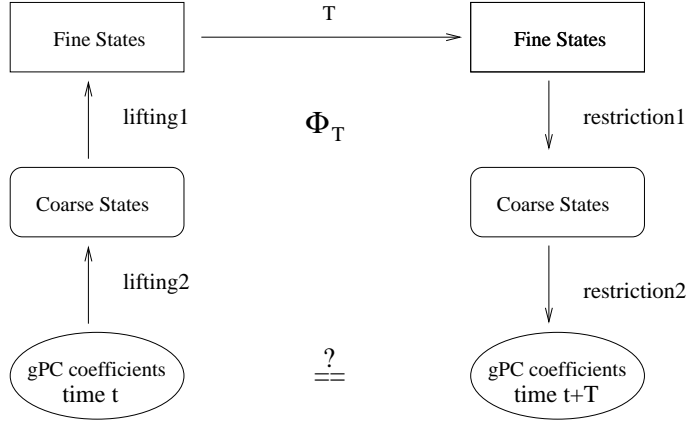,width=.7\textwidth}}
\caption{Fixed-point formulation for coarse random steady states.}
\label{fpmultiscale:fig}
\end{figure}

In equation-free computations of random coarse-scale steady states, we use the
time-stepper for the fully coarse scale states ${\bm X}_{cc}$ (the gPC coefficients),
${\bm \Phi}_T$, to construct a fixed point equation
\begin{equation}
     {\bm X}_{cc} = {\bm \Phi}_T ( {\bm X}_{cc} )
\label{fixedpoint:eqn}
\end{equation}
The operator ${\bm \Phi}_T$ involves repeated lifting and repeated
restriction procedures across two scale gaps, as illustrated in Fig.
\ref{fpmultiscale:fig}.
The solution of Equation (\ref{fixedpoint:eqn}) can be attemped either by direction iteration, by
Newton's method with numerically estimated Jacobians (for a small number of gPC coefficients) or,
more systematically, by matrix-free fixed point algorithms such as Newton-Krylov GMRES
(for large numbers of gPC coefficients) \cite{Kelley:95}.
Once the fully coarse-scale steady states are available,
ensemble realizations and thus probability distributions of coarse-scale steady states can
be immediately obtained by (\ref{PCexpansion:eqn}).
\section{Numerical Results}  \label{results:sec}
Our illustrative example involves the $A+1/2B_2 \rightarrow AB$ reaction,
which can be used as a caricature description of the CO oxidation on a Pt catalyst
surface \cite{Makeev:02}.
Our (intermediate) coarse-scale observables are
the mean coverages of reactants and the vacant catalyst sites.
The fine scale states in each detailed simulation consist of
the numbers of sites occupied by each reactant
and the vacant sites.
In the limit of very large systems, the (intermediate) coarse scale ODEs for mean coverage
based on the particular kinetic mechanism would be
\begin{eqnarray}
   d \theta_A/dt &=& \alpha \theta_* - \gamma \theta_A - k_r \theta_A \theta_B  \cr
   d \theta_B/dt &=& \beta \theta_*^2 - k_r \theta_A \theta_B
\label{coarseODEs:eqn}
\end{eqnarray}
(where $\theta_A$, $\theta_B$ and $\theta_*$ are mean coverages of reactants A, B and vacant sites, respectively).
We will not use these equations in our computations that follow (which are performed
for {\em finite size} systems).

At the fine scale description level, there are four elementary reaction steps:
\begin{eqnarray}
   A(g) +  *_i & \stackrel{\alpha}{\longrightarrow} & A_{*,i}   \cr
   B_2(g) + *_i + *_j & \stackrel{\beta}{\longrightarrow}& B_{*,i} + B_{*,j}  \cr
   A_{*,i} &\stackrel{\gamma}{\longrightarrow}& A(g) + *_i   \cr
   A_{*,i}+B_{*,j} &\stackrel{k_r}{\longrightarrow}& AB(g) + *_i + *_j ;
\end{eqnarray}
here $(g)$ refers to gas phase reactants,
$*_{i(j)}$ a vacant site, and $A(B)_{*,i(j)}$ are the absorbed reactants on the surface.
At the continuum limit, the rates of four reactions are
given respectively by $r_1=\alpha \theta_* N_{tot}$, $r_2=\beta \theta_*^2 N_{tot}$,
$r_3=\gamma \theta_A N_{tot}$ and $r_4=k_r \theta_A \theta_B N_{tot}$,
where $N_{tot}$ is the number of sites on the reacting surface.
For our finite-size system, the four reaction rates at a given time
$t$ (based on which  the reaction probabilities are computed) would
be $r_1 = \alpha N_*$, $r_2 = {1 \over 2} {\beta \over {N_{tot}}}
N_*(N_*-1)$, $r_3 = \gamma N_A$ and $r_4 = {{k_r} \over {N_{tot}}}
N_A N_B$, where $N_A$, $N_B$ and $N_*$ are numbers of sites taken
respectively by $A$, $B$ and vacant slots at time $t$.

Our implementations of this chemical reaction scheme follows the
Gillespie algorithm \cite{Gillespie:77,Gillespie:92}.
The index of reaction that will take place next depends on the random variable $p_1$,
uniformly distributed over the domain $[0,1]$:
reaction $j$ will occur if $\sum_{i=1}^{j-1} r_i / \sum_{i=1}^4 r_i \leq p_1 \leq \sum_{i=1}^j r_i / \sum_{i=1}^4 r_i$.
The time at which this reaction takes place is $-ln(p_2)/\sum_{i=1}^4 r_i$,
where $p_2$ is also a uniformly distributed random variable over $[0,1]$.
For each realization in the fine-level SSA, a stoichiometric matrix
is then used to keep track of the changes in the numbers of the
reactants and vacant sites over time.

\begin{figure}
\centerline{\epsfig{figure=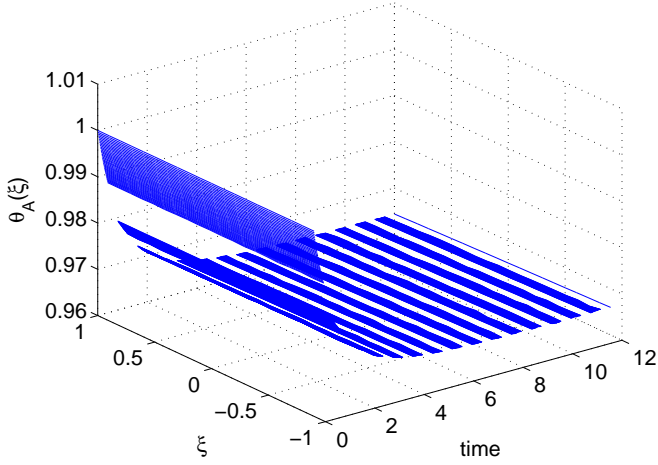,width=.7\textwidth}}
\caption{Evolution
of $\theta_A$ as a function of the random variable $\xi$ through
coarse projective integration} \label{thetaA3d:fig}
\end{figure}

\begin{figure}
\centerline{\epsfig{figure=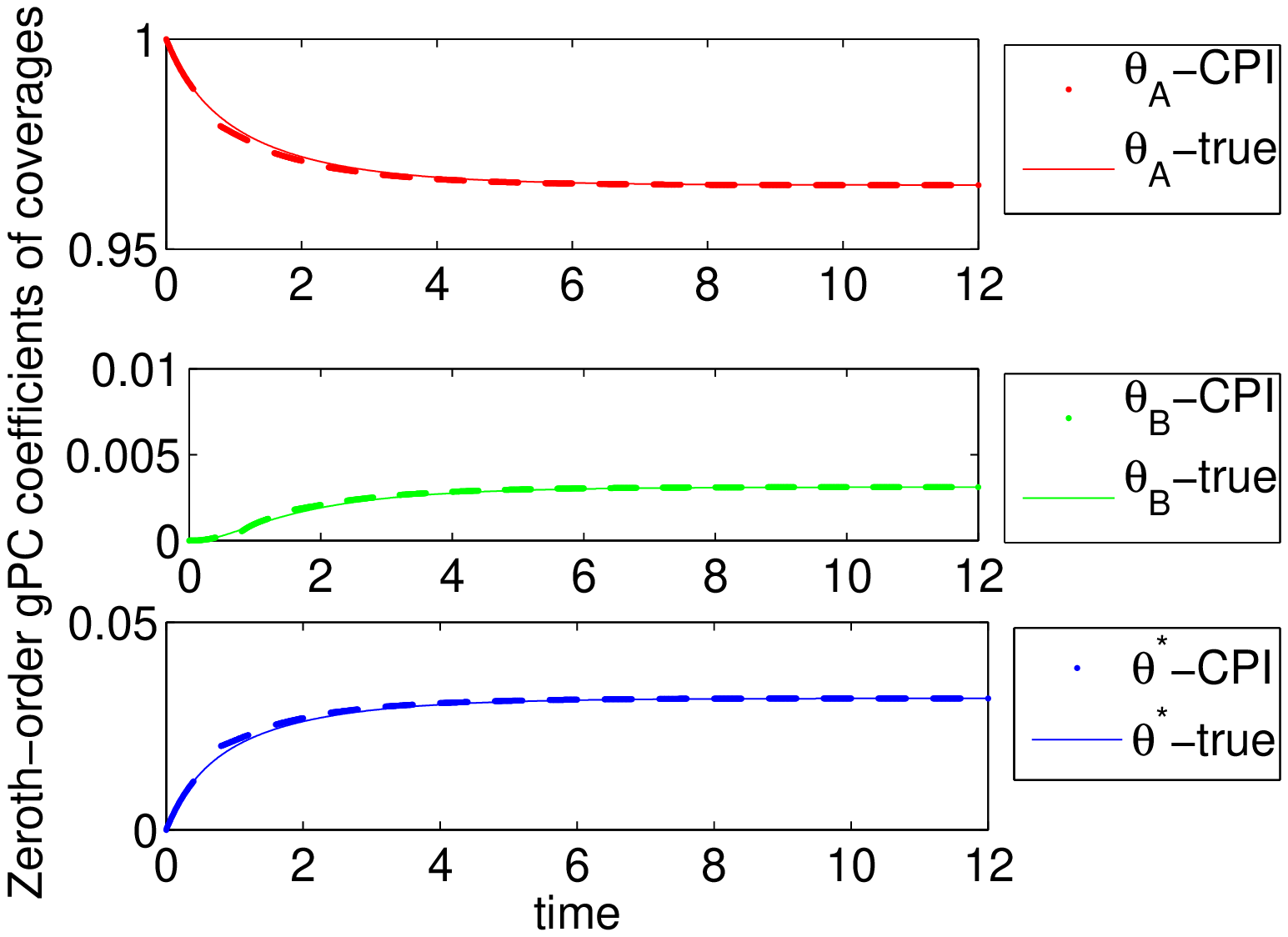,width=.7\textwidth}}
\caption{Evolution of the zeroth-order gPC coefficient of the mean
coverages computed by coarse projective integration, with Ne=40,000
randomly selected $\xi$ values; symbols: CPI results; lines: gPC
coefficients obtained via Monte Carlo simulation of the coarse
ODEs.} \label{cpilsord0:fig}
\end{figure}

\begin{figure}
\centerline{\epsfig{figure=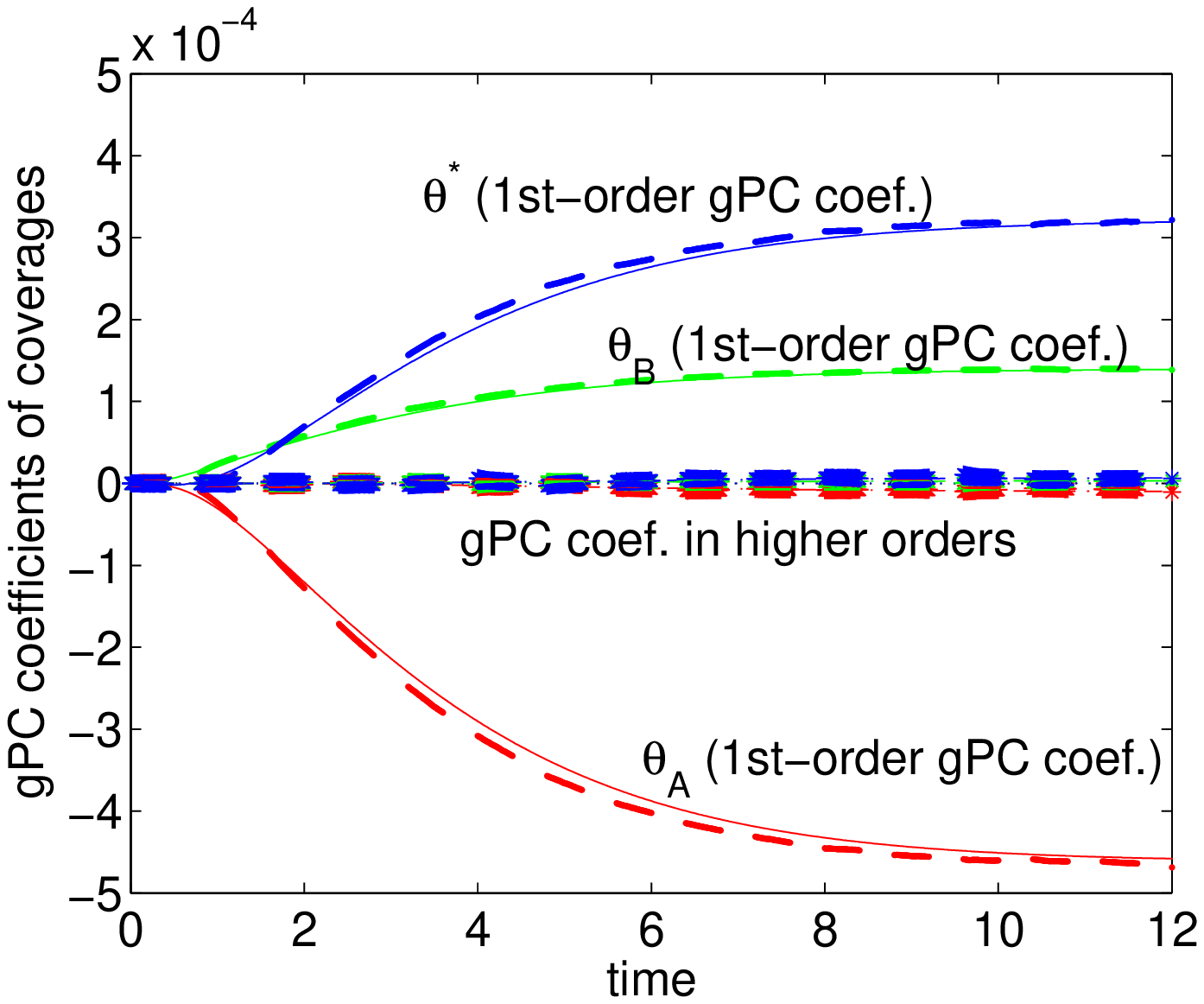,width=.7\textwidth}}
\caption{Evolution of additional (higher-order) gPC coefficients of
the mean coverages computed by coarse projective integration with
Ne=40,000 randomly selected $\xi$ values; symbols: CPI results;
lines: gPC coefficients obtained via Monte Carlo simulation of the
coarse ODEs.} \label{cpilsord123:fig}
\end{figure}

\begin{figure}
\centerline{\epsfig{figure=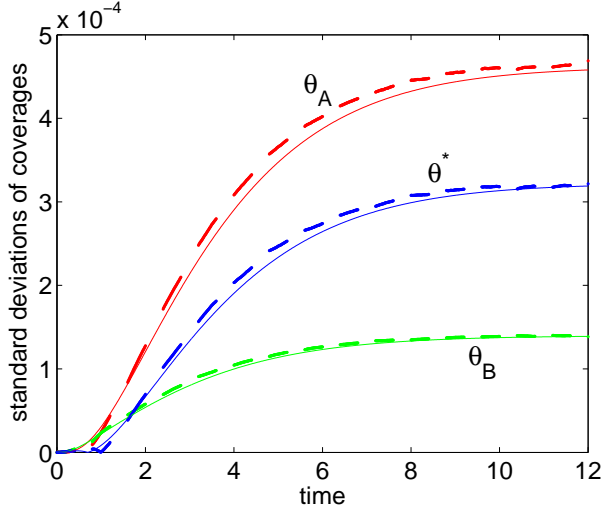,width=.7\textwidth}}
\caption{Evolution of the standard deviations of the mean coverages
computed by coarse projective integration, with Ne=40,000 randomly
selected $\xi$ values; symbols: CPI results; lines: Monte Carlo
simulation of the coarse ODEs.} \label{cpilsstd:fig}
\end{figure}

In the simulations that follow, the parameters $\alpha$, $\gamma$, and $k_r$ are considered known, and
set respectively to $1.6$, $0.04$ and $4$;
the uncertain parameter is $\beta$: $\beta=6.0+0.25\xi$, where $\xi$ is a random
variable uniformly distributed over $[-1,1]$.
Legendre polynomial chaos is chosen as the basis for the random, fully coarse-scale states.
The highest truncation order for these fully coarse observables is chosen as 3 ($P=3$).
A total ensemble of 40,000 realizations are used in lifting and restriction
between the intermediate coarse (mean coverages) to the fully coarse (gPC coefficients) level.
Consistent with each of these 40,000 ensemble realizations,
1000 fine-scale realizations were simulated;
each of them uses $200^2 (N_{tot}=200^2)$ sites on the catalyst surface.
The fine-(intermediate) coarse restriction ${\cal M}$ consists in taking the average
coverage over these 1000 fine-scale realizations corresponding to each (intermediate)
coarse observable (i.e., mean coverage of each reactant).
We perform successive fine-intermediate coarse and then intermediate-fully coarse restrictions
at 40 successive coarse-scale time steps ($\Delta t_c= 0.01s$).
We then use least-squares fitting to estimate temporal derivatives of the
gPC coefficients based on values at the last 5 time steps.
These numerical derivatives are then used (in a simple, forward Euler projective scheme)
to calculate gPC states after a relatively large time step ($\Delta t_{cc} = 0.8s$).
Figure \ref{thetaA3d:fig} shows evolution of the mean coverage of reactant $A$ as a function
of the random variable $\xi$ in the time domain; the ``empty" time intervals in the plot
are the projective forward Euler ``jumps", over which we do not simulate.
Figures \ref{cpilsord0:fig} and \ref{cpilsord123:fig} contrast the evolution of the
fully coarse-scale observables (gPC coefficients), computed by ensemble average
from direct Monte Carlo simulation of the coarse-scale ODEs (\ref{coarseODEs:eqn}),
and also computed through projective integration of the two-scale-gap system;
Fig. \ref{cpilsstd:fig} compares the standard deviations of mean coverages
computed from the gPC coefficients in the two approaches.
The results indicate good agreement between projective integration
computations and true evolutions at the fully coarse gPC level.

\begin{figure}
\centerline{\epsfig{figure=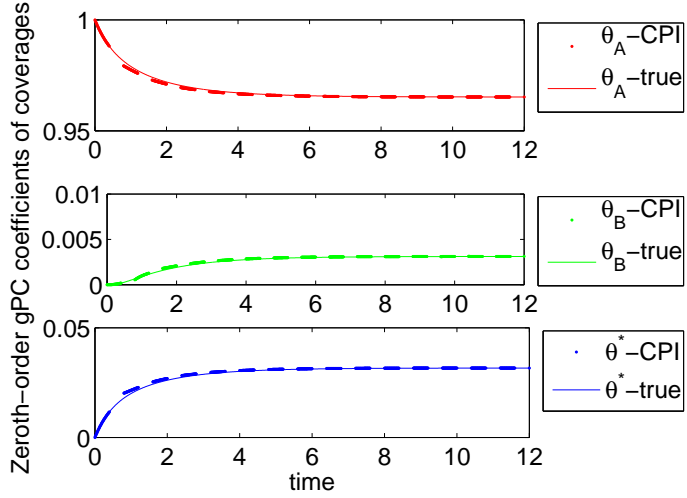,width=.7\textwidth}}
\caption{Evolution of the zeroth-order gPC coefficient of the mean
coverages computed by coarse projective integration for Ne=200 $\xi$
values at Gauss-Legendre points; symbols: CPI results; lines: gPC
coefficients obtained via Monte Carlo simulation of coarse ODEs.}
\label{cpicublsord0:fig}
\end{figure}

\begin{figure}
\centerline{\epsfig{figure=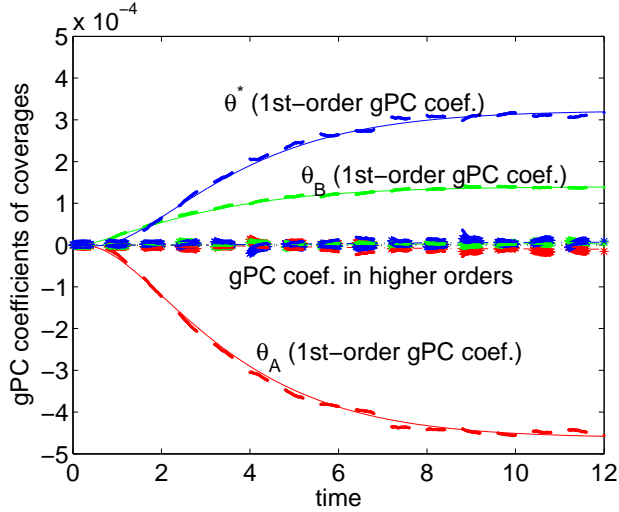,width=.7\textwidth}}
\caption{Evolution of additional (higher-order) gPC coefficients of
the mean coverages computed by coarse projective integration with
Ne=200 $\xi$ values at Gauss-Legendre points; symbols: CPI results;
lines: gPC coefficients obtained via Monte Carlo simulation of the
coarse ODEs.} \label{cpicublsord123:fig}
\end{figure}

\begin{figure}
\centerline{\epsfig{figure=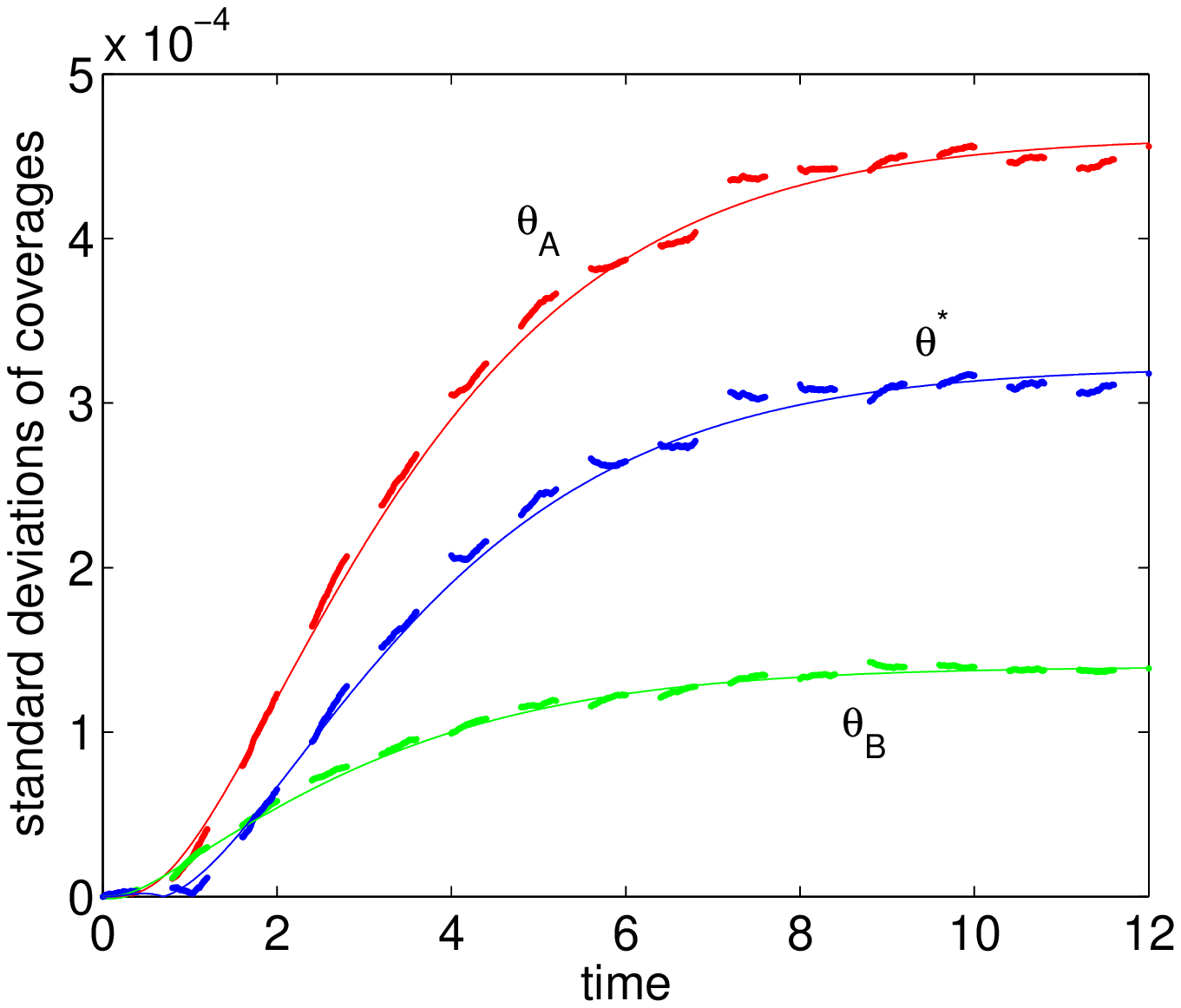,width=.7\textwidth}}
\caption{Evolution of the standard deviations of the mean coverages
computed by coarse projective integration with Ne=200 $\xi$ values
at Gauss-Legendre points; symbols: CPI results; lines: Monte Carlo
simulation of coarse ODEs.} \label{cpicublsstd:fig}
\end{figure}

\begin{figure}
\centerline{\epsfig{figure=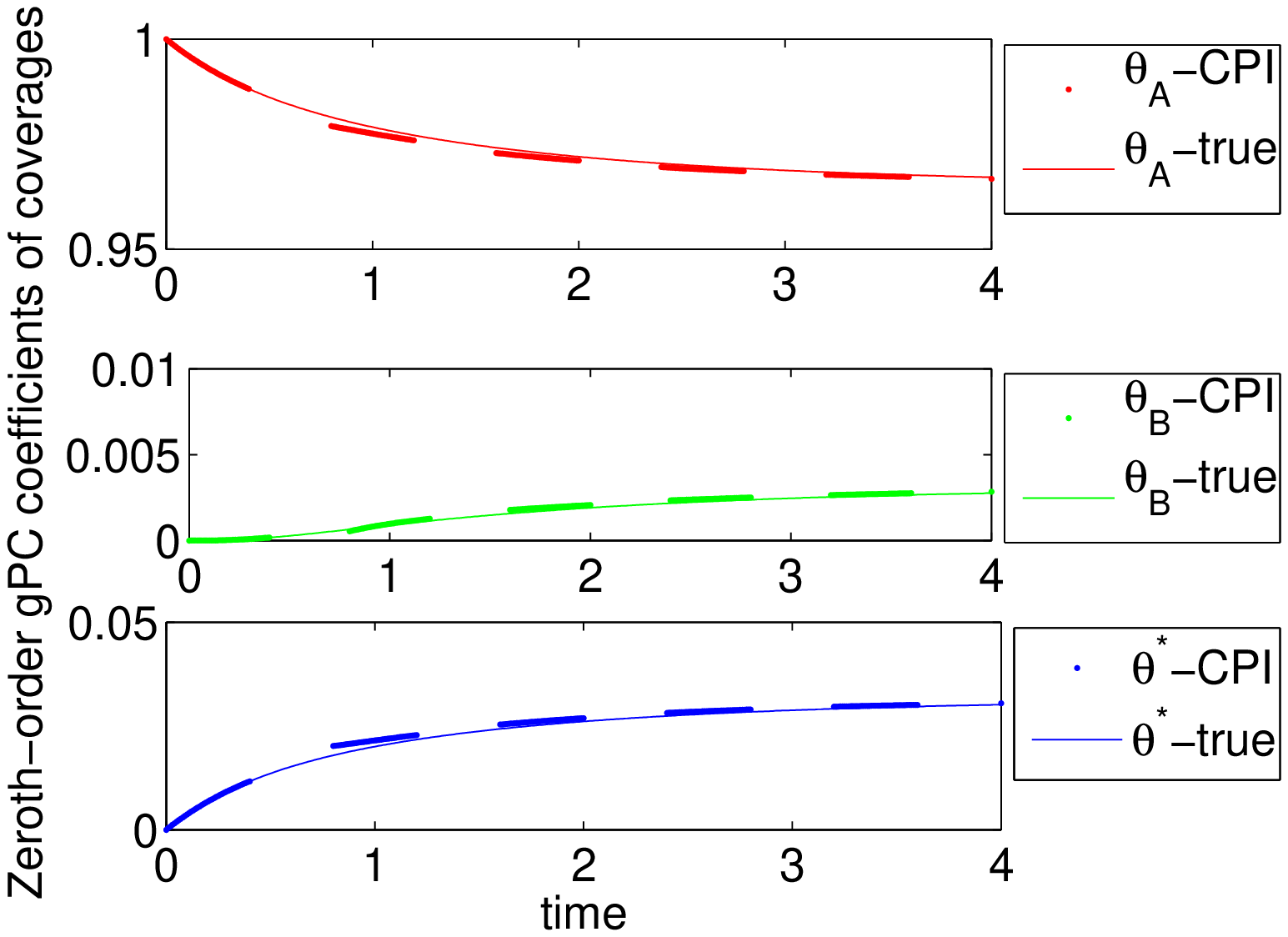,width=.7\textwidth}}
\caption{Evolution of the zeroth-order gPC coefficients of the mean
coverages computed by coarse projective integration for Ne=200 {\em
randomly selected} $\xi$ values; symbols: CPI results; lines: gPC
coefficients obtained via Monte Carlo simulation of coarse ODEs with
Ne=40,000.} \label{cpilsord0Ne200:fig}
\end{figure}

\begin{figure}
\centerline{\epsfig{figure=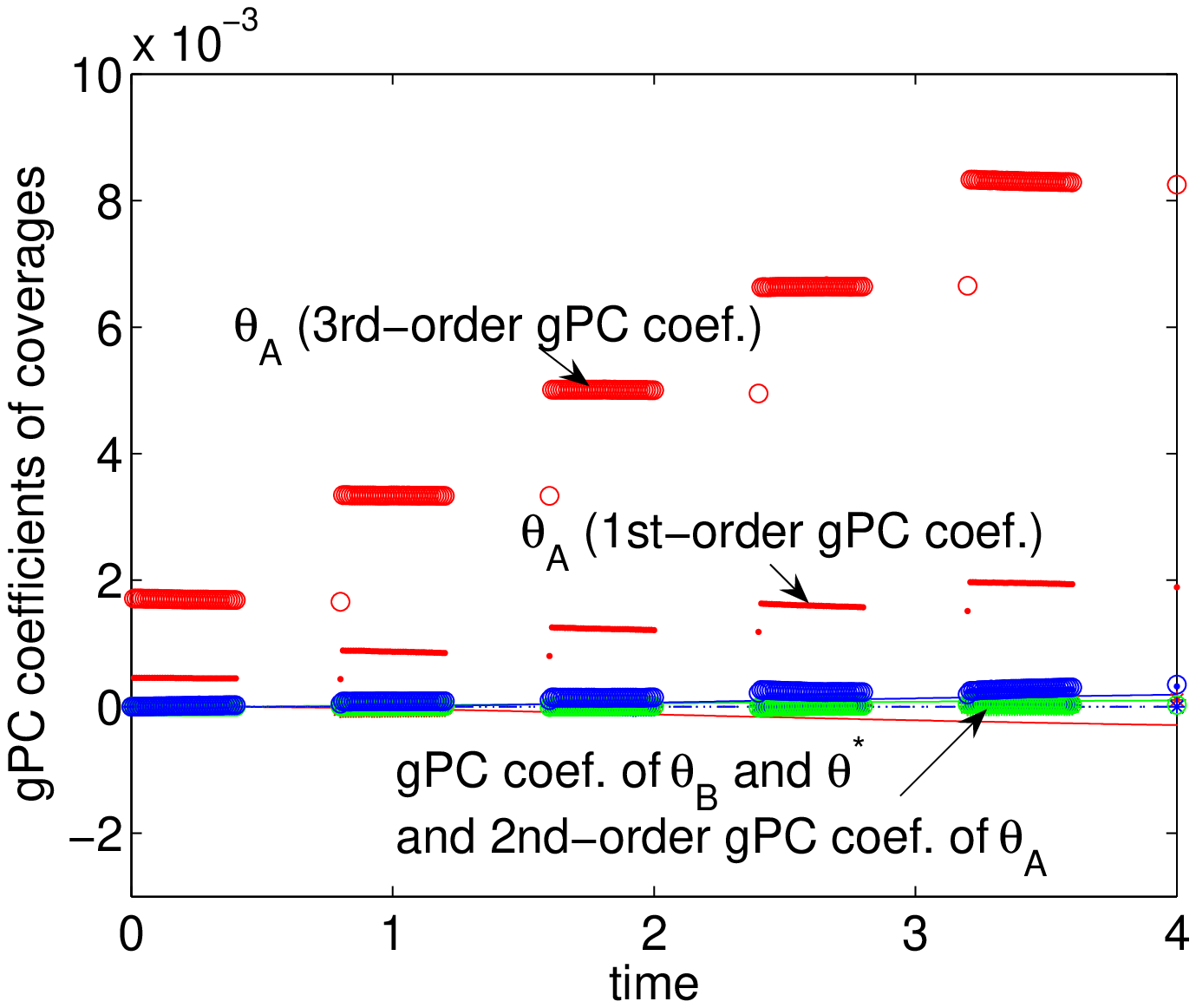,width=.7\textwidth}}
\caption{Evolution of additional, higher-order gPC coefficients of
the mean coverages computed by coarse projective integration with
Ne=200 {\em randomly selected} $\xi$ values. Red circles and dots represent the 
3rd- and 
1st-order gPC coefficients of $\theta_A$ obtained by CPI, respectively. Jumps during the 
course of their evolution are caused by insufficient number of coarse 
realizations of coverages used to perform ensemble average when 
computing the gPC coefficiens. 
Clustered green and blue symbols in the lower part of the figure stand for gPC 
coefficients of reactant $\theta_B$ and vacancy
$\theta^*$, respectively. 
Lines represent gPC coefficients obtained via Monte Carlo simulation of the
coarse ODEs with Ne=40,000. Red, green and blue solid lines stand for
respectively the 1st-order gPC coefficients of $\theta_A$, $\theta_B$ and 
$\theta^*$. Red, green and blue dot-dashed lines stand for
respectively the 2nd-order gPC coefficients of $\theta_A$, $\theta_B$ and 
$\theta^*$. Red, green and blue dotted lines stand for
respectively the 3rd-order gPC coefficients of $\theta_A$, $\theta_B$ and 
$\theta^*$.} \label{cpilsord123Ne200:fig}
\end{figure}

In an attempt to further accelerate the computation, the Gauss-Legendre quadrature was used to
approximate the inner product in the intermediate coarse to fully coarse restriction (\ref{PCcoefficient:eqn}).
In the corresponding fully-to-intermediate coarse lifting, only coarse states corresponding to values
of $\xi$ {\em at the Gauss-Legendre points} were generated.
Figures \ref{cpicublsord0:fig}-\ref{cpicublsstd:fig} show coarse projective integration results
using this method, which utilizes only 200 coarse-scale realizations (rather than 40,000 ones).
To display the effectiveness of this technique, we implemented the original Monte-Carlo simulations
in the intermediate level by using a total of 200 realizations; the results are shown in figures \ref{cpilsord0Ne200:fig}
and \ref{cpilsord123Ne200:fig}.
With such a small number of realizations in a Monte Carlo simulation,
while the zeroth-order gPC coefficients of mean coverages of reactants are well captured,
other, higher order coefficients deviate significantly from their true trajectories even at the
beginning of the projective integration.
This is clearly due to the lack of sufficient realizations of mean coverage
when ensemble averaging is used to compute the corresponding gPC coefficients.
Lifting only around Gauss-Legendre quadrature points can enhance the effectiveness of our
approach to simulate the evolution of gPC coefficients by using a much smaller number of intermediate coarse-level realizations
than that required by the standard Monte-Carlo simulation.
Computing the gPC coefficient evolution by sampling only close to
qudarature points has been proposed (for Gauss-Hermite quadrature,
in a problem of approximating PC coefficients for random temperature
distribution) in \cite{Maitre:02}.

\begin{figure*}
\centerline{\epsfig{figure=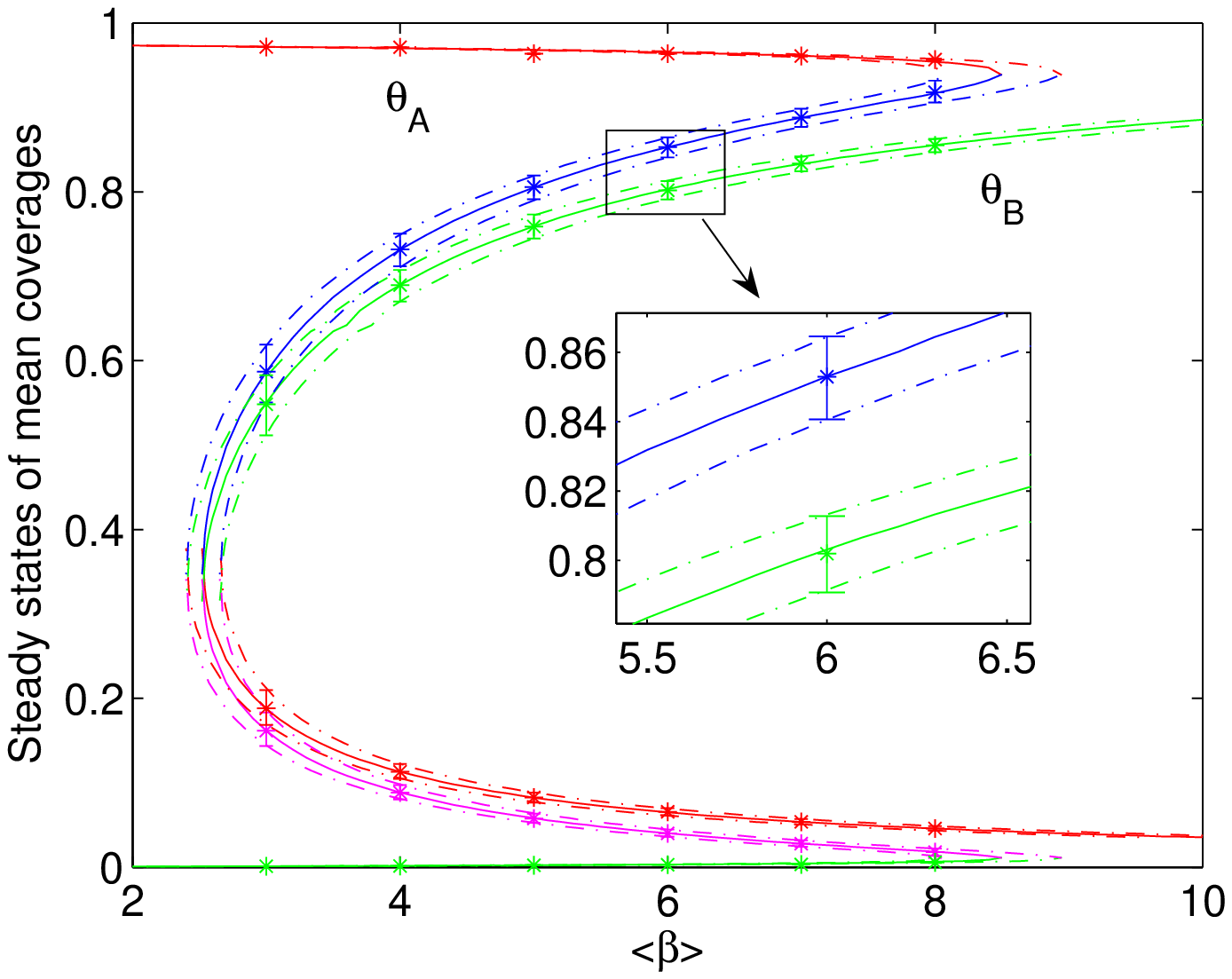,width=.7\textwidth}}
\caption{Bifurcation diagram of the random model steady states with
respect to the mean $<\beta>$ of the uncertain parameter
distribution. Coarse fixed point computation and continuation with
Ne=40,000 randomly selected $\xi$ values; red and green objects:
stable steady states; blue and magenta objects: unstable steady
states. See text for the inset.} \label{randomfp:fig}
\end{figure*}

\begin{figure*}
\centerline{\epsfig{figure=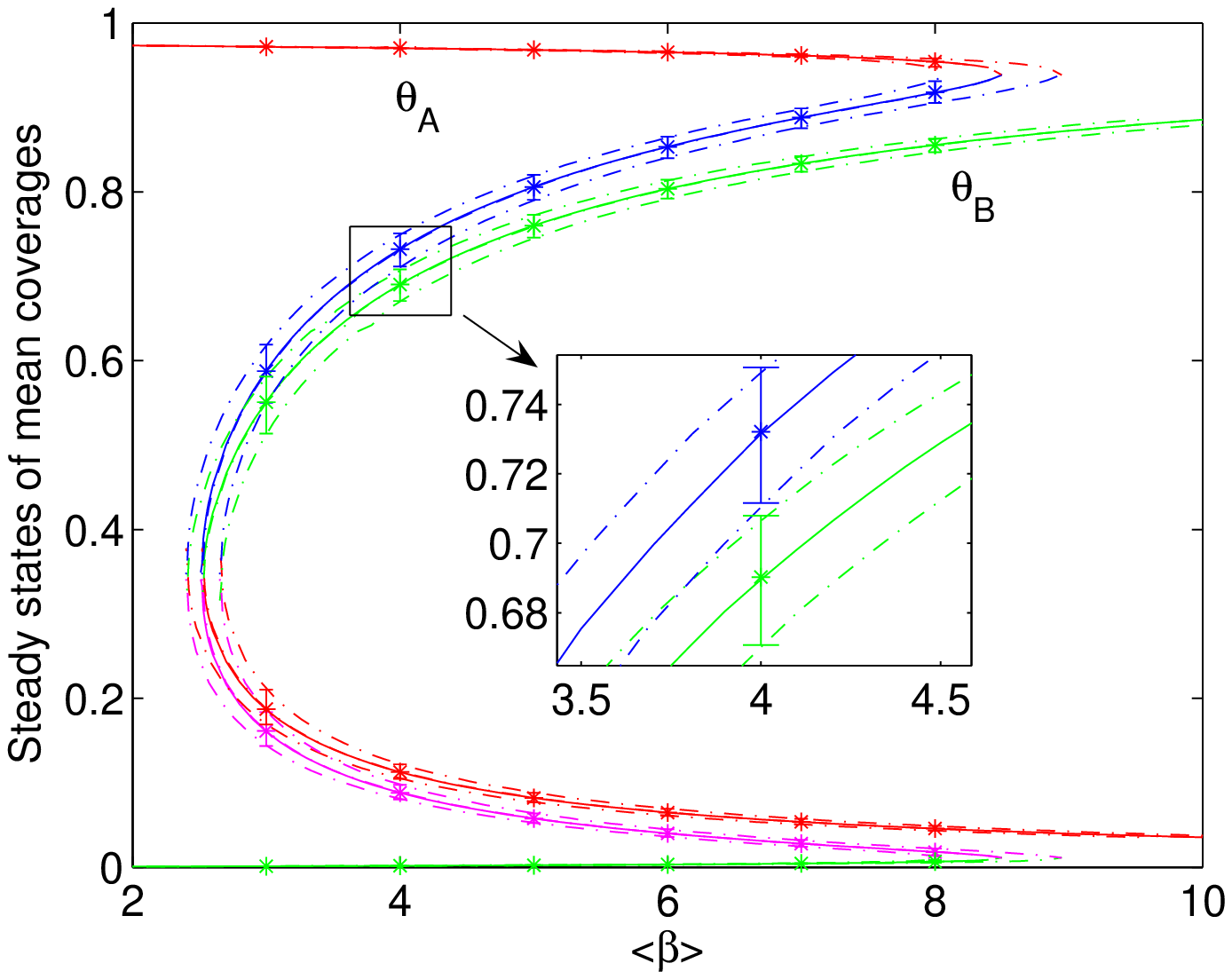,width=.7\textwidth}}
\caption{Bifurcation diagram of the random model steady states with
respect to the mean $<\beta>$ of the uncertain parameter
distribution. Coarse fixed point computation and continuation with
Ne=200 $\xi$ values at Gauss-Legendre points; red and green objects:
stable steady states; blue and magenta objects: unstable steady
states. See text for the inset.} \label{randomfpgauss:fig}
\end{figure*}

We also use a matrix-free Newton-Krylov GMRES method to
converge on the (deterministic) stable/unstable fully coarse steady states of the gPC coefficient description,
out of which random (intermediate) coarse steady state distributions can be obtained by (\ref{PCexpansion:eqn}).
In this computation, $\beta=<\beta>(1.0+0.05\xi)$, where $\xi$ is again a uniform random
variable in $[-1,1]$.
In Fig. \ref{randomfp:fig}, dashed lines are true envelopes of random coarse-scale steady states
(computed by setting $\beta=1.05<\beta>$ and $\beta=0.95<\beta>$),
and solid lines are the derterministic steady states when $\beta=<\beta>$.
Error bars and stars represent respectively ranges and means of random steady states of
mean coverages computed using the matrix-free, time-stepper based Newton-Krylov-GMRES method.
Clearly, this equation-free fixed-point computation can correctly
reproduce the random steady states.
Again, values of $\xi$ at Gauss-Legendre points are used to generate realizations of intermediately-coarse-scale
observables (mean coverages) and implement fixed-point computations (Fig. \ref{randomfpgauss:fig}).
The random steady states can be accurately captured by this technique as well, while the computational load
decreases significantly.
The approach has been linked with a continuation algorithm to trace the bifurcation diagrams of
the random steady states as a function of $<\beta>$; observe that {\em unstable}
random steady states can thus be computed, and bifurcation points (such as turning points
of random steady states) in parameter space can be approximated (see Figures \ref{randomfp:fig} and \ref{randomfpgauss:fig}).

\section{Summary} \label{conclusion:sec}
Equation-free methods were used at two successive layers, in this
paper, to enable uncertainty quantification computations on models
of reacting systems for which no coarse-grained, continuum
description is available in closed form.
Coarse projective integration was used to accelerate the computation
of transient dynamics of the problem solution distributions (at gPC
coefficients level).
Random stable/unstable steady state computation, and their
parametric/bifurcation analysis at this gPC coefficient level was
also demonstrated.
Gauss quadrature rules were used to effectively reduce the computational load while
preserving the accuracy of results of the equation-free results.
The method should be applicable in UQ for a wider class of reacting
problems, beyond those described through a fine-scale SSA simulator,
for which no closed-form kinetic equations are available.
We believe that the approach can still serve in cases where we know
how to describe uncertainty in the {\em microscopic} simulation
parameters, but we cannot easily translate that in uncertainty
descriptions for parameters at the coarse-grained level.

{\bf Acknowledgements}. This work was partially supported by DARPA and by the US DOE. I.G.K. also
acknowledges the support of a 2005 Guggenheim Fellowship.


\bibliographystyle{unsrt}
\bibliography{UQGLPMMS}

\end{document}